\documentclass[11pt]{amsart}
\usepackage{amsmath}
\usepackage{amssymb}
\usepackage{amsthm}
\usepackage{mathrsfs}
\usepackage{enumerate}
\usepackage[english]{babel}
\usepackage{hyperref}
\usepackage[all,cmtip]{xy}
\entrymodifiers={+!!<0pt,\fontdimen22\textfont2>}

\newtheorem{thm}{Theorem}[section]
\newtheorem{lem}[thm]{Lemma}
\newtheorem*{thm-i}{Theorem}
\newtheorem{cor}[thm]{Corollary}
\newtheorem{prop}[thm]{Proposition}

\theoremstyle{definition}

\newtheorem{rem}[thm]{Remark}
\newtheorem{defn}[thm]{Definition}

\newtheorem*{acknowledgments*}{Acknowledgments}

\numberwithin{equation}{section}

\theoremstyle{remark}

\mathchardef\ordinarycolon\mathcode`\: 
\def\vcentcolon{\mathrel{\mathop\ordinarycolon}} 
\providecommand*\coloneqq{\mathrel{\vcentcolon\mkern-1.2mu}=}

\def\Z{{\mathbb Z}} 

\DeclareMathOperator\res{res}

\def\KK{{K\!K}}
\def\K{{\widetilde{K}}}
\def\E{{\widetilde E}}

\def\L{{\widetilde E}}

\def\F{{\mathcal F}}
\def\H{{\mathcal H}}

\def\I{{\mathcal I}}
\def\C{{\mathcal C}}
\def\O{{\mathcal O}}

\def\cast{$C^{*}$}
\def\pt{*}

\newcommand\thmA{(2)}
\newcommand\thmC{(1)}

\begin{document}

\title{Restriction maps in equivariant $\KK$-theory.}
\author{Otgonbayar Uuye}
\date{December 13, 2010}
\address{
Department of Mathematical Sciences\\
University of Copenhagen\\
Universitetsparken 5\\
DK-2100 Copenhagen E\\
Denmark}
\email{otogo@math.ku.dk}
\urladdr{http://www.math.ku.dk/~otogo}
\keywords{$K$-theory, $\KK$-theory}
\subjclass[2010]{Primary (19K35); Secondary (46L80)}
\maketitle

\begin{abstract} We extend McClure's results on the restriction maps in equivariant $K$-theory to bivariant $K$-theory:

Let $G$ be a compact Lie group and $A$ and $B$ be $G$-\cast-algebras. Suppose that $\KK^{H}_{n}(A, B)$ is a finitely generated $R(G)$-module for every $H \le G$ closed and $n \in \Z$. Then, if $\KK^{F}_{*}(A, B) = 0$ for all $F \le G$ {\em finite cyclic}, then $\KK^{G}_{*}(A, B) = 0$.
\end{abstract}

\setcounter{section}{-1}
\section{Introduction}

Let $G$ be a compact Lie group and let $X$ be a finite $G$-CW-complex. For any closed subgroup $H \le G$, we have a restriction functor in equivariant $K$-theory: 
	\begin{equation}
	\res^{G}_{H}: K_{G}(X) \to K_{H}(X).		
	\end{equation}

As an application of the generalized Atiyah-Segal completion theorem of \cite{MR935523}, McClure proved the following.
\begin{thm-i}[{McClure \cite[Theorem A]{MR862427}}] If $x \in K_{G}(X)$ restricts to zero in $K_{H}(X)$ for every {\em finite} subgroup $H$ of $G$, then $x = 0$. 
\end{thm-i}

Combining with Jackowski's result \cite[Corollary 4.3]{MR0448377}, one obtains the following.
\begin{thm-i}[{Jackowski-McClure \cite[Corollary C]{MR862427}}] 
If $K_{F}^{*}(X) = 0$ for all $F \le G$ finite cyclic, then $K_{G}^{*}(X) = 0$.
\end{thm-i}

We extend these to bivariant $K$-theory as follows. 

\begin{thm}\label{thm KK} Let $G$ be a compact Lie group and $A$ and $B$ be $G$-\cast-algebras. Suppose that $\KK^{H}_{n}(A, B)$ is a finitely generated $R(G)$-module for every $H \le G$ closed and $n \in \Z$.
\begin{enumerate}
\item[\thmC] 
Then, if $\KK^{F}_{*}(A, B) = 0$ for all $F \le G$ {\em finite cyclic}, then $\KK^{G}_{*}(A, B) = 0$.
\item[\thmA] Suppose, in addition, that $\KK^{F}_{n}(A, B)$ is a finitely generated group for all $F \le G$ finite and $n \in \Z$. Then, if $x \in \KK^{G}(A, B)$ restricts to zero in $\KK^{H}(A, B)$ for all $H \le G$ {\em finite}, then $x = 0$. 
\end{enumerate}
\end{thm}

\begin{rem}
\begin{enumerate}
\item See \cite{MR2083579} for a dual result for restriction maps in $K$-homology of spaces with actions of discrete groups.
\item Theorem \ref{thm KK} is in stark contrast to the results of Heath Emerson, where he showed that even for circle actions, noncommutative algebras can behave very differently from commutative ones. \cite{2010arXiv1004.2970E}
\end{enumerate}
\end{rem}

In fact, we prove the following. This is done mainly for clarity, but as an added bonus, we see that Theorem~\ref{thm KK} holds for equivariant $E$-theory as well.
\begin{thm}\label{thm L}
Let $G$ be a compact Lie group and let $\L_{G}^{*}$ be an $RO(G)$-gradable module theory over $\K_{G}^{*}$. Suppose that $\L_{H}^{n}(S^{0})$ is a finitely generated $R(G)$-module for every $H \le G$ closed and $n \in \Z$. Let $X$ be a finite based $G$-CW-complex.

\begin{enumerate}
\item[\thmC] Then, if $\L_{F}^{*}(X) = 0$ for all $F \le G$ {\em finite cyclic}, then $\L_{G}^{*}(X) = 0$.

\item[\thmA] Suppose, in addition, that $\L_{F}^{n}(S^{0})$ is a finitely generated group for all $F \le G$ finite and $n \in \Z$. Then, if $x \in \L_{G}^{*}(X)$ restricts to zero in $\L_{H}^{*}(X)$ for all  $H \le G$ {\em finite}, then $x = 0$. 
\end{enumerate} 
\end{thm}

The proof follows \cite{MR862427} very closely. In Section \ref{sec RO(G)}, we show that Theorem~\ref{thm L} implies Theorem~\ref{thm KK}. In Section \ref{sec Completion}, we extend the generalized Atiyah-Segal completion theorem of \cite{MR935523}, supplying the missing ingredient needed to finish the proof in Section \ref{sec Proof of main thm}. 

\begin{acknowledgments*} The author thanks the Centre for Symmetry and Deformation at the University of Copenhagen and the Danish National Research Foundation for support.
\end{acknowledgments*}

\section{$RO(G)$-graded cohomology theories}\label{sec RO(G)}

Let $G$ be a compact Lie group. A based $G$-space is a $G$-space with a $G$-fixed base point. In the rest of the paper, we assume that all $G$-spaces are $G$-CW-complexes and all cohomology theories are equivariant and reduced cohomology theories.

For a finite-dimensional representation $V$ of $G$, we write $S^{V}$ for the one-point compactification of $V$, considered a based $G$-space with base point the point at infinity.

\subsection{$RO(G;U)$-gradable theories}
We fix a complete universe $U$. (cf.\ \cite[Definition IX.2.1]{MR1413302}). 

\begin{defn} An {\em $RO(G)$-graded} cohomology theory is an $RO(G; U)$-graded cohomology theory in the sense of \cite[Definition XIII.1.1]{MR1413302}. A {\em $\Z$-graded} cohomology theory is an $RO(G; U^{G})$-graded cohomology theory (any trivial universe would work). We say that a $\Z$-graded cohomology theory is {\em $RO(G)$-gradable} if it is the $\Z$-graded part of an $RO(G)$-graded theory.
\end{defn}


Let $\E_{G}^{*}$ be a $\Z$-graded cohomology theory. 
For a closed subgroup $H \le G$ and a based $H$-CW-complex $X$, we define 
	\begin{equation}
	\E_{H}^{*}(X) \coloneqq \E_{G}^{*}(G_{+} \wedge_{H}X).	
	\end{equation}
Then $\E_{H}^{*}$ is a $\Z$-graded cohomology theory on based $H$-spaces. If $X$ is actually a based $G$-CW-complex, then we have a natural $G$-equivariant  identification 
	\begin{equation}
	G_{+} \wedge_{H}X \cong G/H_{+} \wedge X
	\end{equation}
and the collapse map $G/H \to \pt$ gives rise to a natural transformation
	\begin{equation}
	\res^{G}_{H}: \E_{G}^{*} \to \E_{H}^{*}
	\end{equation}
called the {\em restriction} map. 

\subsection{Bivariant $K$-theory} The following is the main example we have in mind. First note that $\K_{G}^{*}$ is an $RO(G)$-graded commutative ring theory with $\K_{G}^{V}(X) = \KK^{G}(C_{0}(S^{V}), C_{0}(X))$ and $R(G) \cong \K_{G}(S^{0})$.

\begin{prop} Let $G$ be a compact Lie group and let $A$ and $B$ be $G$-\cast-algebras. For a finite based $G$-CW-complex $X$ and finite-dimensional real representation $V$ of $G$, we define
	\begin{equation}
	\L_{G}^{V}(X) \coloneqq \KK^{G}(A \otimes C_{0}(S^{V}), B \otimes C_{0}(X)).
	\end{equation}
Then the following holds.
\begin{enumerate}[(i)]
\item $\L_{G}^{*}$ defines an $RO(G)$-graded cohomology theory on the category of finite based $G$-CW-complexes.
\item $\L_{G}^{*}$ extends to an $RO(G)$-graded cohomology theory on the category of based $G$-CW-complexes.
\item $\L_{G}^{*}$ is a module theory over $\K_{G}^{*}$.
\end{enumerate}
\end{prop}
\begin{proof}
(i) See \cite{MR918241}.
(ii) By Adams' representation theorem \cite[Theorem XIII.3.4]{MR1413302}, $\L_{G}^{*}$ is represented by an $\Omega$-$G$-prespectrum, hence extends to an $RO(G)$-graded cohomology theory on the category of $G$-CW-complexes. See \cite{MR1193150}.
(iii) The module structure
	\begin{equation}
	\L_{G}^{V}(X) \times \K^{W}_{G}(Y) \to \L_{G}^{V+W}(X \wedge Y).
	\end{equation}
is given by the Kasparov product 
	\begin{align}
	\KK^{G}(A(S^{V}), B(X)) \times \KK^{G}(C_{0}(S^{W}), C_{0}(Y))\\ \to \KK^{G}(A(S^{V+W}), B(X \wedge Y)).
	\end{align}
\end{proof}

It is well-known that for $H \le G$,  
	\begin{equation}
	\KK^{G}(A, B \otimes C_{0}(G/H_{+})) \cong \KK^{H}(A, B)	
	\end{equation}
and the restriction map is induced by $G/H_{+} \to S^{0}$. Hence we obtain the following corollary.

\begin{cor}
Suppose that Theorem~\ref{thm L} holds. Then Theorem~\ref{thm KK} holds.
\qed
\end{cor}
\section{Atiyah-Segal Completion}\label{sec Completion}

First we abstract the main finiteness condition from Theorem~\ref{thm L}.
\begin{defn} Let $R$ be a unital commutative ring and let $\L_{G}^{*}$ be a $\Z$-graded cohomology theory with values in $R$-modules. We say that $\L_{G}^{*}$ is {\em finite} over $R$ if $\L_{G}^{n}(X)$ is a finitely generated $R$-module for every {\em finite} based $G$-CW-complex $X$ and $n \in \Z$.
\end{defn}

Clearly, this is equivalent to asking that $\L_{H}^{k-n}(S^{0}) \cong \L_{G}^{k}(G/H_{+} \wedge S^{n})$ is a finitely generated $R$-module for $H \le G$.

\begin{lem}\label{lem cohomology} Let $G$ be a compact Lie group and let $R$ be a unital commutative ring. Let $\L_{G}^{*}$ be a $\Z$-graded cohomology theory with values in $R$-modules. Suppose that $R$ is Noetherian and $\L_{G}^{*}$ is finite over $R$. Then for any family $\I$ of ideals in $R$, the following defines a $\Z$-graded cohomology theory with values in pro-$R$-modules:
	\begin{equation}
	\L_{G}^{*}(X)^{\wedge}_{\I} \coloneqq \{\L_{G}^{*}(Y)/J \cdot \L_{G}^{*}(Y)\}.
	\end{equation}
where $Y \subseteq X$ runs over the finite based $G$-CW-subcomplexes of $X$ and $J$ runs over the finite products of ideals in $\I$.  
\end{lem}
Note that in this lemma, it is enough to have $\L_{G}^{*}$ to be a cohomology theory on {\em finite} based $G$-CW-complexes (only finite wedges are considered in the additivity axiom). 
\begin{proof} Exactness follows from the Artin-Rees lemma. See the proof of \cite[Lemma 2.1]{MR935524}.
\end{proof}

\subsection{Bott Periodicity}
Let $V$ be a complex $G$-representation. By Bott periodicity \cite[Theorem 4.3]{MR0228000}, $\K_{G}^{0}(S^{V})$ is a free $\K_{G}^{0}(S^{0})$-module generated by the Bott element $\lambda_{V} \in \K_{G}^{0}(S^{V})$. The {\em Euler class} of $V$ is defined to be $e^{*}(\lambda_{V}) \in \K_{G}^{0}(S^{0})$, where $e: S^{0} \to S^{V}$ is the obvious map. 

\begin{lem}\label{lem Bott} Let $\L_{G}^{*}$ be an $RO(G)$-graded module theory over $\K_{G}^{*}$. Then for any complex representation $V$, multiplication by the Bott element $\lambda_{V} \in \K_{G}^{0}(S^{V})$ gives an isomorphism 
	\begin{equation}
	\L_{G}^{0}(S^{0}) \cong \L_{G}^{0}(S^{V}). 
	\end{equation}
If $V \subseteq W$ are complex representations and $i: S^{V} \to S^{W}$ is the inclusion, then the following diagram commutes
	\begin{equation}
	\xymatrix{
	\L_{G}^{0}(S^{0}) \ar[r] \ar[d]^{\cdot \chi_{W-V}}& \L_{G}^{0}(S^{W}) \ar[d]^{i^{*}}\\
	\L_{G}^{0}(S^{0}) \ar[r] & \L_{G}^{0}(S^{V})\\
	}.
	\end{equation}
\end{lem}
\begin{proof} Let $\lambda_{V}^{-1} \in \K_{G}^{V}(S^{0})$ denote the inverse Bott element: it has the property that
	\begin{equation}
	\lambda_{V} \cdot \lambda_{V}^{-1} = \lambda_{V}^{-1} \cdot \lambda_{V} = 1 \in \K_{G}^{V}(S^{V}) \cong \K_{G}^{0}(S^{0}).
	\end{equation}
Then multiplication by $\lambda_{V}^{-1}$ gives the inverse map
	\begin{equation}
	\L_{G}^{0}(S^{V}) \to \L_{G}^{V}(S^{V}) \cong \L_{G}^{0}(S^{0}).
	\end{equation}
The second statement is shown for $\L_{G}^{*} = \K_{G}^{*}$ in \cite[page 4]{MR935523}. The general case follows by functoriality.
\end{proof}
	
\subsection{Completion}
A class of subgroups of $G$ closed under subconjugacy is called a {\em family}. A family $\C$ of subgroups of $G$ determines a class, again denoted  $\C$, of ideals of $R(G)$ by the kernels of the restriction maps:
	\begin{equation}
	\ker(\res^{G}_{H}: R(G) \to R(H)),\quad H \in \C,	
	\end{equation}
hence a topology on any $R(G)$-module. 

The following is a straightforward generalization of \cite[Theorem 3.1]{MR935523}.
\begin{thm}\label{thm AS} Let $G$ be a compact Lie group and let $\L_{G}^{*}$ be an $RO(G)$-gradable module theory over $\K_{G}^{*}$, which is finite over $R(G)$. 

Let $\C$ be a family of subgroups of $G$. For any based $G$-CW-complex $X$, if $\L_{H}^{*}(X)^{\wedge}_{\C|H} = 0$ for all $H \in \C$, then $\L_{G}^{*}(X)^{\wedge}_{\C} =0.$

\end{thm}
\begin{proof} By \cite[Corollary 3.3]{MR0248277}, $R(G) = \K_{G}^{0}(S^{0})$ is Noetherian. Hence, by Lemma~\ref{lem cohomology}, $\L_{G}^{*}(X)^{\wedge}_{\C}$ is a cohomology theory.

Now the proof of \cite[Theorem 3.1]{MR935523} carries over ad verbatum, once we extend Bott periodicity to $\L_{G}^{*}$ as in  Lemma~\ref{lem Bott}.
\end{proof}

\begin{cor}\label{cor AS}
Let $E\C$ denote the classifying space of $\C$. For any finite based $G$-CW-complex $X$, the projection map $E\C_{+} \to S^{0}$ gives completion
	\begin{equation}
	\L_{G}^{*}(E\C_{+} \wedge X) \cong \lim \L_{G}^{*}(Y \wedge X) \cong \lim \L_{G}^{*}(X)^{\wedge}_{\C},
	\end{equation}
where $Y$ runs over finite based subcomplexes of $E\C_{+}$. 	
\end{cor}
\begin{proof}
The inverse system $\L_{G}^{*}(X)^{\wedge}_{\C}$ satisfies the Mittag-Leffler condition and $\L_{G}^{*}(Y \wedge X)$ is $\C$-complete for any finite based subcomplex $Y \subset E\C_{+}$ (cf.\ \cite[Corollary 2.1]{MR935523}).
\end{proof}

\section{Proof of Theorem \ref{thm L}}\label{sec Proof of main thm}
\subsection{$\F$-spaces}

Let $\F$ be a family of subgroups of $G$. We say that a based $G$-CW-complex $X$ is an $\F$-space if all the isotropy groups, except at the base point, are in $\F$. The following lemma says that in the proof of Theorem~\ref{thm L}, we may assume that $X$ is an $\F$-space, for any $\F$ containing all finite cyclic subgroups of $G$.

\begin{lem}\label{lem F-space} 
Let $G$ be a compact Lie group and let $\L_{G}^{*}$ be an $RO(G)$-gradable module theory over $\K_{G}^{*}$, which is finite over $R(G)$. 

Let $\F$ be a family containing all finite cyclic subgroups of $G$. Then for any finite based $G$-CW-complex $X$,  the top horizontal map in the commutative diagram
	\begin{equation}
	\xymatrix{\L_{G}^{*}(X) \ar[r] \ar[d] & \lim_{Y \subset E\F_{+}}\L_{G}^{*}(Y \wedge X) \ar[d]\\
	\prod_{F \in \F} \L_{F}^{*}(X) \ar[r] &\lim_{Y \subset E\F_{+}}\prod_{F \in \F} \L_{F}^{*}(Y \wedge X)\\ 
	},
	\end{equation}
is injective. Here $Y$ runs over the finite based subcomplexes of $E\F_{+}$, the horizontal maps are induced by the projections $Y \wedge X \to X$ and the vertical maps are restrictions. 
\end{lem}
\begin{proof} The $\F$-topology on $\L_{G}^{*}(X)$ is Hausdorff by \cite[Corollary 3.3]{MR862427}. Hence, the claim follows from Corollary~\ref{cor AS}.
\end{proof}

Let $\C$ denote the family of finite cyclic subgroups of $G$.

\begin{proof}[Proof of Theorem~\ref{thm L}\thmC]
By assumption, $\L_{F}^{*}(X) = 0$ for all $F \in \C$. Let $Y$ be a finite based $G$-CW-complex, which is a $\C$-space. Then the zero skeleton $Y^{0}$ and the skeletal quotients $Y^{n}/Y^{n-1}$ are finite wedges of $G$-spaces of the form $G/F_{+} \wedge S^{n}$ with $F \in \C$. It follows that $\L_{G}^{*}(Y \wedge X) = 0$. Hence by Lemma~\ref{lem F-space}, $\L_{G}^{*}(X) = 0$.
\end{proof}

\subsection{Induction}
We write $\O_{G}$ for the category whose objects are orbit spaces $G/H$, where $H \le G$ is a closed subgroup, and whose morphisms are homotopy classes of $G$-maps. 

Recall that a compact Lie group is said to {\em cyclic} if it has a topological generator (an element whose powers are dense) and {\em hyperelementary} if it is an extension of a cyclic group by a finite $p$-group.

We write $\H$ for the class of hyperelementary subgroups of $G$ and let $\O_{\H}$ denote the full subcategory of $\O_{G}$ of orbits $G/H$ with $H$ subconjugate to a subgroup in $\H$.

\begin{lem}\label{lem ind} Let $G$ be a compact Lie group and let $\L_{G}^{*}$ be an $RO(G)$-gradable module theory over $\K_{G}^{*}$. Then, for any based $G$-CW-complex,  the restriction maps induce an isomorphism
	\begin{equation}
	\L_{G}^{*}(X) \cong \lim_{\O_{\H}} \L_{H}^{*}(X).
	\end{equation}
\end{lem}
\begin{proof} Follows from Propositions 2.1 and 2.2 of \cite{MR862427}.
\end{proof}

For any abelian group $M$, let $M^{\wedge}_{\Z}$ denote its adic completion $\lim_{n} M/nM$.

\begin{proof}[Proof of Theorem \ref{thm L}\thmA]
Let $\F$ denote the family of finite subgroups of $G$.

By Lemma \ref{lem ind}, we may assume that $G$ is a hyperelementary group and by Lemma \ref{lem F-space}, we may assume that $X$ an $\F$-space. 

Let $G$ be a hyperelementary group and $X$ an $\F$-space. Then the restriction map 
	\begin{equation}
	\L_{G}^{*}(X)^{\wedge}_{\Z} \to \lim_{F \in \O_{\F}} \L_{F}^{*}(X)^{\wedge}_{\Z}.
	\end{equation}
is an isomorphism by \cite[Theorem 1.1]{MR862427}. By \cite[Corollary 3.3]{MR862427}, the adic topologies on $\L_{G}^{*}(X)$ and $\L_{F}^{*}(X)$ are Hausdorff. This completes the proof.
\end{proof}
\bibliographystyle{amsalpha}
\bibliography{../BibTeX/biblio}
\end{document}